\documentclass[12pt]{amsart}
\usepackage{amsmath,amssymb,amstext,amsthm,amscd,epic,eepic}

\theoremstyle{plain}
\newtheorem{thm}{Theorem}
\newtheorem*{thm2}{Theorem}

\theoremstyle{definition}
\newtheorem*{ack}{Acknowledgements}
\theoremstyle{remark}
\newtheorem*{rem}{Remark}
\newtheorem*{rems}{Remarks}

\def\({\left(}
\def\){\right)}

\def\fun{\rightarrow}
\def\lfun{{\longrightarrow}}
\def\inc{\hookrightarrow}

\def\l0{\lambda^0}
\def\li{\lambda^\infty}
\def\psii{\psi_\infty}

\def\bC{{\mathbb{C}}}
\def\bP{{\mathbb{P}}}

\def\cM{{\mathcal M}}
\def\ocM{{\overline{\cM}}}

\begin{document}

\title{Recursive formula for
$\psi^g-\lambda_1\psi^{g-1}+\cdots+(-1)^g\lambda_g$ in $\ocM_{g,1}$}

\author{D.\ Arcara}
\address{Department of Mathematics, University of Utah,
155 S. 1400 E., Room 233, Salt Lake City, UT 84112-0090, USA}
\email{arcara@math.utah.edu}

\author{F.\ Sato}
\address{School of Mathematics, Korean Institute for Advanced Study,
Cheongnyangni 2-dong, Dongdaemun-gu, Seoul 130-722, South Korea}
\email{fumi@kias.re.kr}

\begin{abstract}
Mumford proved that
$\psi^g-\lambda_1\psi^{g-1}+\cdots+(-1)^g\lambda_g=0$ in the Chow ring of
$\cM_{g,1}$ [Mum83].
We find an explicit recursive formula for
$\psi^g-\lambda_1\psi^{g-1}+\cdots+(-1)^g\lambda_g$
in the tautological ring of $\ocM_{g,1}$ as a combination of classes supported
on boundary strata.
\end{abstract}

\maketitle

\section{Introduction}

Mumford proved in [Mum83] that
$\psi^g-\lambda_1\psi^{g-1}+\cdots+(-1)^g\lambda_g=0$ in the Chow ring of
$\cM_{g,1}$.
Moreover, he showed that this class is supported on the boundary strata with a
marked genus $0$ component.
Graber and Vakil proved in [GraVak05] that every codimension $g$ class in the
tautological ring of $\ocM_{g,1}$ is supported on the boundary strata with at
least one genus $0$ component.

We complement these results by finding an explicit recursive formula for
$\psi^g-\lambda_1\psi^{g-1}+\cdots+(-1)^g\lambda_g$ in the tautological ring
of $\ocM_{g,1}$ as a combination of classes supported on boundary strata.
It is clear from the formula being recursive that all the boundary strata
have a genus $0$ component in them, but it is not obvious from the formula
that the marked point must be on a genus $0$ component.
We simplified the formula for $g<5$ in Section \ref{simplify}, and checked
that this is the case.

\begin{thm}\label{thm}
In the tautological ring of $\ocM_{g,1}$,
$$ \sum_{i=0}^g (-1)^i \lambda_i \psi^{g - i} =
\sum_{h=1}^g \( 1 - \frac{h}{g} \) \iota_h{}_*(c_h), $$
where
$$ c_h := \sum_{i=0}^{g-1} (-1)^{h+i} \left[
\( \sum_{j=0}^h (-1)^j \l0_j \psi_0^{i-j} \)
\( \sum_{j=0}^{g-h} (-1)^j \li_j \psii^{g-1-i-j} \) \right], $$
$\iota_h$ is the natural boundary map
$$ \iota_h \colon \ocM_{h,2} \times \ocM_{g-h,1} \lfun \ocM_{g,1}, $$
$\psi_0$, $\psii$ are descendents at the marked points glued by $\iota_h$, and
$\l0$, $\li$ are the $\lambda$-classes on $\ocM_{h,2}$ and $\ocM_{g-h,1}$,
respectively.
\end{thm}

This formula is actually the first step of an algorithm which calculates each
of the classes $\psi^g$, $\lambda_1\psi^{g-1}$, $\dots$, $\lambda_g$ in terms
of classes supported on boundary strata.
We want to single out the class
$\psi^g-\lambda_1\psi^{g-1}+\cdots+(-1)^g\lambda_g$, though, because it is the
only class we found so far in the tautological ring of $\ocM_{g,1}$ which has
a nice recursive formula, and can therefore be easily calculated.

\begin{ack}
We would like to thank A.\ Bertram for many valuable discussions.
This paper would have never been possible without his explanations,
suggestions and guidance throughout the project.
\end{ack}

\section{Virtual localization}

The main tool we use to prove our theorems is the virtual localization theorem
by Graber and Pandharipande [GraPan99].

\begin{thm2}[Virtual localization theorem]
Suppose $f\colon X\fun X'$ is a $\bC^*$-equivariant map of proper
Deligne--Mumford quotient stacks with a $\bC^*$-equivariant perfect
obstruction theory.
If $i'\colon F'\inc X'$ is a fixed substack and $c\in A^*_{\bC^*}(X)$, let
$f|_{F_i}\colon F_i\fun F'$ be the restriction of $f$ to each of the fixed
substacks $F_i\subseteq f^{-1}(F')$.
Then 
$$\sum_{F_i}{f|_{F_i}}_*\frac{i_{F_i}^*c}{\epsilon_{\bC^*}(F_i^{\text{vir}})}
=\frac{i'^*f_*c}{\epsilon_{\bC^*}(F'^{\text{vir}})}, $$
where $i_{F_i}\colon F_i\fun X$ and $\epsilon_{\bC^*}(F^{\text{vir}})$ is the
virtual equivariant Euler class of the ``virtual'' normal bundle
$F^{\text{vir}}$.
\end{thm2}

\begin{rem}
The conditions in the theorem are satisfied for the Kontsevich--Manin spaces
$\ocM_{g,n}(\bP^m,d)$ of stable maps, and $\epsilon_{\bC^*}(F^{\text{vir}})$
can be explicitly computed in terms of $\psi$ and $\lambda$-classes [GraPan99]
(see also [FabPan05]).
\end{rem}

We define a $\bC^*$-action on $\bP^{1}$ by
$ a \cdot [x : y] = [x : a y] $
for $a\in\bC^*$ and $[x:y]\in\bP^{1}$.
There are two fixed points, $0$ and $\infty$, and the torus acts with weight
$1$ on the tangent space at $0$ and $-1$ on the tangent space at $\infty$.
This $\bC^*$-action induces $\bC^*$-actions on $\ocM_{g,n}(\bP^{1},d)$, and we
shall consider the trivial $\bC^*$-action on $\ocM_{g,n}$.

\section{Proof of Theorem \ref{thm}}

We use virtual localization on the natural function
$ f \colon \ocM_{g,3}(\bP^1,1) \lfun \ocM_{g,3} \times (\bP^1)^3 $
defined by
$ f([g \colon (C,p_1,p_2,p_3) \fun \bP^1]) = 
((C_{\textrm{stab}},p_1,p_2,p_3),g(p_1),g(p_2),g(p_3)). $
Consider the fixed locus
$$ F' := \ocM_{g,4} \times \{ 0 \} \times \{ \infty \} \times \{ \infty \}
\inc \ocM_{g,3} \times (\bP^1)^3, $$
and apply the virtual localization theorem with $c=[1]^\textrm{vir}$ to obtain
$$ \sum_{F_i} \( f|_{F_i} \)_*
\frac{[1]^\textrm{vir}}{\epsilon_{\bC^*}(F_i^{\text{vir}})} =
\frac{i'^* f_*[1]^\textrm{vir}}{t(-t)(-t)}. $$

There are $g+1$ fixed loci mapping to $F'$.
One fixed locus has a marked point mapping to $0$ and a curve in $\ocM_{g,3}$
mapping to $\infty$.
We shall denote it by $F_0$.
Then there are $g$ fixed loci which have a curve in $\ocM_{h,2}$ mapping to
$0$ and a curve in $\ocM_{g-h,3}$ mapping to $\infty$ (with $1\leq h\leq g$).
We shall denote these fixed loci by $F_h$.
Note that $F_0\simeq\ocM_{g,3}$ and $F_h\simeq\ocM_{h,2}\times\ocM_{g-h,3}$
($1\leq h\leq g$).

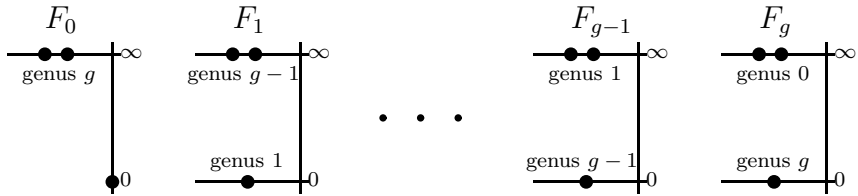
\begin{eqnarray*}
\setlength{\unitlength}{0.01cm}
\begin{picture}(1120,250)(0,0)
\thicklines
\put(50,220){$F_0$}
\put(140,0){\line(0,1){200}}
\put(140,15){\circle*{15}}
\put(0,185){\line(1,0){150}}
\put(50,185){\circle*{15}}
\put(80,185){\circle*{15}}
{\tiny
\put(150,10){$0$}
\put(150,180){$\infty$}
\put(20,150){$\textrm{genus }g$}
}
\put(300,220){$F_1$}
\put(390,0){\line(0,1){200}}
\put(250,15){\line(1,0){150}}
\put(320,15){\circle*{15}}
\put(250,185){\line(1,0){150}}
\put(300,185){\circle*{15}}
\put(330,185){\circle*{15}}
{\tiny
\put(400,10){$0$}
\put(400,180){$\infty$}
\put(270,35){$\textrm{genus }1$}
\put(240,150){$\textrm{genus }g-1$}
}
\put(500,100){\circle*{7}}
\put(600,100){\circle*{7}}
\put(550,100){\circle*{7}}
\put(750,220){$F_{g-1}$}
\put(840,0){\line(0,1){200}}
\put(700,15){\line(1,0){150}}
\put(770,15){\circle*{15}}
\put(700,185){\line(1,0){150}}
\put(750,185){\circle*{15}}
\put(780,185){\circle*{15}}
{\tiny
\put(850,10){$0$}
\put(850,180){$\infty$}
\put(720,150){$\textrm{genus }1$}
\put(690,35){$\textrm{genus }g-1$}
}
\put(1000,220){$F_g$}
\put(1090,0){\line(0,1){200}}
\put(950,15){\line(1,0){150}}
\put(1020,15){\circle*{15}}
\put(950,185){\line(1,0){150}}
\put(1000,185){\circle*{15}}
\put(1030,185){\circle*{15}}
{\tiny
\put(1100,10){$0$}
\put(1100,180){$\infty$}
\put(970,35){$\textrm{genus }g$}
\put(970,150){$\textrm{genus }0$}
}
\end{picture}
\end{eqnarray*}

Since $i'^*f_*[1]^\textrm{vir}$ is a polynomial in $t$, the sum of the
contributions from the coefficient of $t^{-4}$ on each fixed locus is $0$.
Call this contribution $a_{-4}$.
We have that $\pi_{1,*}(\pi_{2,*}(a_{-4}\cdot\psi_3))=0$.
We now calculate the contribution to the left hand side one fixed locus at the
time.

\begin{itemize}
\item
For $F_0$, we obtain
$$ \frac{[1]^\textrm{vir}}{\epsilon_{\bC^*}(F_0^{\text{vir}})} = \frac{1}{t}
\cdot (-1)^g
\frac{t^g + \lambda_1 t^{g-1} + \cdots + \lambda_g}{- t (- t - \psii)}, $$
and the coefficient of $t^{-4}$ is
$ - \( \psii^{g+1} - \lambda_1 \psii^g + \cdots +
(-1)^g \lambda_g \psii \). $
Under the isomorphism $F_0\simeq\ocM_{g,3}$, $\psii$ gets identified with
$\psi_1$, and the contribution is therefore
$$ - \( \psi_1^{g+1} - \lambda_1 \psi_1^g + \cdots +
(-1)^g \lambda_g \psi_1 \). $$
\item
For $F_h$ ($1\leq h\leq g$), we obtain
$$ \frac{[1]^\textrm{vir}}{\epsilon_{\bC^*}(F_h^{\text{vir}})} =
\frac{t^h - \l0_1 t^{h-1} + \cdots + (-1)^h \l0_h}{t (t - \psi_0)} \cdot
(-1)^{g-h}
\frac{t^{g-h} + \li_1 t^{g-h-1} + \cdots + \li_{g-h}}{- t (- t - \psii)}, $$
and the coefficient of $t^{-4}$ is\footnote{Note that $c'_h$ is the
summation (with the appropriate sign) of all possible products of codimension
$g$ of a class on the curve mapping to $0$ with a class on the curve mapping
to $\infty$.}
$$ c'_h := \sum_{i=0}^g (-1)^{h+i} \left[
\( \sum_{j=0}^h (-1)^j \l0_j \psi_0^{i-j} \)
\( \sum_{j=0}^{g-h} (-1)^j \li_j \psii^{g-i-j} \) \right]. $$
This is a class of codimension $g$ in $\ocM_{h,2}\times\ocM_{g-h,3}$ which
maps to the codimension $g+1$ class $\iota_h{}_*(c'_h)$ in $\ocM_{g,3}$ under
$(f|_{F_h})_*$.
\end{itemize}

To summarize, we obtain that
$$ - \( \psi_1^{g+1} - \lambda_1 \psi_1^g + \cdots +
(-1)^g \lambda_g \psi_1 \) + \sum_{h=1}^g \iota_h{}_*(c'_h) = 0 $$
in $\ocM_{g,3}$.
The first step is now to multiply by $\psi_3$ and push-forward to
$\ocM_{g,2}$.

\begin{itemize}
\item
If $h=0$, we obtain
$ - 2 g \( \psi_1^{g+1} - \lambda_1 \psi_1^g + \cdots +
(-1)^g \lambda_g \psi_1 \) $
in $\ocM_{g,2}$.
\item
If $1\leq h<g$, note that, since the third marked point is on the curve at
$\infty$, we are really multiplying by $\psi_3$ in $\ocM_{g-h,3}$ and
pushing-forward to $\ocM_{g-h,2}$.
We therefore obtain, by Dilaton, the class
$ 2 (g - h) \iota_h{}_*(c'_h), $
which is a class of codimension $g+1$ in $\ocM_{g,2}$.
\item
If $h=g$, then $\psi_3=0$ because it is a descendent at a marked point of a
genus $0$ curve with $3$ markings (the curve mapping to $\infty$).
\end{itemize}

Let us now suppose that $h<g$.
The second and last step is to push-forward this class via the map that
forgets the second marked point.

\begin{itemize}
\item
If $h=0$, we obtain, by String,
$ - 2 g \( \psi_1^g - \lambda_1 \psi_1^{g - 1} + \cdots +
(-1)^g \lambda_g \). $
\item
If $1\leq h<g$, we obtain, by String, the class
$ 2 (g - h) \iota_h{}_*(c_h), $
where $c_h$ is just $c'_h$ with every power of $\psii$ lowered by $1$ (with
the convention that $\psii^{-1}=0$), i.e.,
$$ c_h = \sum_{i=0}^{g-1} (-1)^{h+i} \left[
\( \sum_{j=0}^h (-1)^j \l0_j \psi_0^{i-j} \)
\( \sum_{j=0}^{g-h} (-1)^j \li_j \psii^{g-1-i-j} \) \right]. $$
\end{itemize}

Putting it all together, we obtain that
$$ - 2 g \( \psi_1^g - \lambda_1 \psi_1^{g - 1} + \cdots + (-1)^g \lambda_g \)
+ \sum_{h=1}^g 2 (g - h) \iota_h{}_*(c_h) = 0, $$
from which we can derive the formula of Theorem \ref{thm}.
\qed

\begin{rems}
(I) By taking the coefficient of $t^{-3-j}$ with $j>1$, it is possible to
find a similar formula for
$ \psi^{g+j-1} - \lambda_1 \psi^{g+j-2} + \cdots +
(-1)^g \lambda_g \psi^{j-1} $
in terms of classes supported on boundary strata.

(II) In [GraVak05], Graber and Vakil proved that a codimension $g$ class in
the tautological ring of $\ocM_{g,1}$ can be written as a sum of classes
supported on boundary strata with at least one genus $0$ component.
By induction on $g$, it is easy to see that this is the case for our $c_h$
classes.

(III) Using the same function $f$ as above, but with the fixed locus
$\ocM_{g,2}\times\{0\}^2\times\{\infty\}$ instead of
$\ocM_{g,2}\times\{0\}\times\{\infty\}^2$, it is possible to obtain the
following tautological relation on $\ocM_{g,1}$:
$$ \sum_{h=1}^{g-1} (2h) \iota_h{}_* (c_h) + (2g) \pi_* \( \psi_2^{g+1} -
\lambda_1 \psi_2^g + \cdots + (-1)^g \lambda_g \psi_2 \) = 0. $$
\end{rems}

\section{Explicit formulas for low genus}\label{simplify}

The formula of Theorem \ref{thm} can be simplified recursively, and we
calculated the answer for low values of $g$.
Note that these formulas were already known for $g=1$ and $g=2$, but they were
unknown for higher $g$'s.

Genus $1$: In $\ocM_{1,1}$,
\vspace{-.5cm}
{\tiny
\begin{eqnarray*}
\setlength{\unitlength}{0.0075cm}
\begin{picture}(280,140)(0,0)
\thicklines
\put(100,50){\ellipse{200}{100}}
\path(50,54)(60,52)(70,50)(80,48)(90,46)(100,45)(110,46)(120,48)(130,50)(140,52)(150,54)
\path(70,50)(80,52)(90,54)(100,55)(110,54)(120,52)(130,50)
\put(100,22){\circle*{10}}
\put(50,115){$\psi-\lambda_1$}
\put(210,45){$=0.$}
\end{picture}
\end{eqnarray*}
}

Genus $2$: In $\ocM_{2,1}$,
\vspace{-.5cm}
{\tiny
\begin{eqnarray*}
\setlength{\unitlength}{0.0075cm}
\begin{picture}(900,140)(0,0)
\thicklines
\path(0,50)(1,57)(2,60)(5,66)(8,70)(14,76)(18,79)(21,81)(23,82)(26,84)(34,88)(37,89)(40,90)(52,94)(60,96)(65,97)(72,98)(80,99)(100,100)(101,100)
\path(0,50)(1,43)(2,40)(5,34)(8,30)(14,24)(18,21)(21,19)(23,18)(26,16)(34,12)(37,11)(40,10)(52,6)(60,4)(65,3)(72,2)(80,1)(100,0)(101,0)
\path(175,88)(160,90)(148,94)(140,96)(135,97)(128,98)(120,99)(100,100)(99,100)
\path(175,12)(160,10)(148,6)(140,4)(135,3)(128,2)(120,1)(100,0)(99,0)
\path(350,50)(349,57)(348,60)(345,66)(342,70)(336,76)(332,79)(329,81)(327,82)(324,84)(316,88)(313,89)(310,90)(298,94)(290,96)(285,97)(278,98)(270,99)(250,100)(249,100)
\path(350,50)(349,43)(348,40)(345,34)(342,30)(336,24)(332,21)(329,19)(327,18)(324,16)(316,12)(313,11)(310,10)(298,6)(290,4)(285,3)(278,2)(270,1)(250,0)(249,0)
\path(175,88)(190,90)(202,94)(210,96)(215,97)(222,98)(230,99)(250,100)(251,100)
\path(175,12)(190,10)(202,6)(210,4)(215,3)(222,2)(230,1)(250,0)(251,0)
\path(50,54)(60,52)(70,50)(80,48)(90,46)(100,45)(110,46)(120,48)(130,50)(140,52)(150,54)
\path(70,50)(80,52)(90,54)(100,55)(110,54)(120,52)(130,50)
\path(300,54)(290,52)(280,50)(270,48)(260,46)(250,45)(240,46)(230,48)(220,50)(210,52)(200,54)
\path(280,50)(270,52)(260,54)(250,55)(240,54)(230,52)(220,50)
\put(175,50){\circle*{10}}
\put(50,115){$\psi^2-\lambda_1\psi+\lambda_2$}
\put(360,45){$=$}
\put(500,50){\ellipse{200}{100}}
\path(450,54)(460,52)(470,50)(480,48)(490,46)(500,45)(510,46)(520,48)(530,50)(540,52)(550,54)
\path(470,50)(480,52)(490,54)(500,55)(510,54)(520,52)(530,50)
\put(650,50){\circle{100}}
\path(700,50)(690,48)(680,46)(670,44)(660,42)(650,41)(640,42)(630,44)(620,46)(610,48)(600,50)
\path(690,52)(680,54)\path(670,56)(660,58)\path(650,59)(640,58)\path(630,56)(620,54)\path(610,52)(600,50)
\put(650,20){\circle*{10}}
\put(800,50){\ellipse{200}{100}}
\path(750,54)(760,52)(770,50)(780,48)(790,46)(800,45)(810,46)(820,48)(830,50)(840,52)(850,54)
\path(770,50)(780,52)(790,54)(800,55)(810,54)(820,52)(830,50)
\end{picture}
\end{eqnarray*}
}

Genus $3$: In $\ocM_{3,1}$,
{\tiny
\begin{eqnarray*}
\setlength{\unitlength}{0.0075cm}
\begin{picture}(1750,200)(0,-50)
\thicklines
\path(0,50)(1,57)(2,60)(5,66)(8,70)(14,76)(18,79)(21,81)(23,82)(26,84)(34,88)(37,89)(40,90)(52,94)(60,96)(65,97)(72,98)(80,99)(100,100)(101,100)
\path(0,50)(1,43)(2,40)(5,34)(8,30)(14,24)(18,21)(21,19)(23,18)(26,16)(34,12)(37,11)(40,10)(52,6)(60,4)(65,3)(72,2)(80,1)(100,0)(101,0)
\path(175,88)(160,90)(148,94)(140,96)(135,97)(128,98)(120,99)(100,100)(99,100)
\path(175,12)(160,10)(148,6)(140,4)(135,3)(128,2)(120,1)(100,0)(99,0)
\path(175,88)(190,90)(202,94)(210,96)(215,97)(222,98)(230,99)(250,100)(251,100)
\path(175,12)(190,10)(202,6)(210,4)(215,3)(222,2)(230,1)(250,0)(251,0)
\path(325,88)(310,90)(298,94)(290,96)(285,97)(278,98)(270,99)(250,100)(249,100)
\path(325,12)(310,10)(298,6)(290,4)(285,3)(278,2)(270,1)(250,0)(249,0)
\path(325,88)(340,90)(352,94)(360,96)(365,97)(372,98)(380,99)(400,100)(401,100)
\path(325,12)(340,10)(352,6)(360,4)(365,3)(372,2)(380,1)(400,0)(401,0)
\path(500,50)(499,57)(498,60)(495,66)(492,70)(486,76)(482,79)(479,81)(477,82)(474,84)(466,88)(463,89)(460,90)(448,94)(440,96)(435,97)(428,98)(420,99)(400,100)(399,100)
\path(500,50)(499,43)(498,40)(495,34)(492,30)(486,24)(482,21)(479,19)(477,18)(474,16)(466,12)(463,11)(460,10)(448,6)(440,4)(435,3)(428,2)(420,1)(400,0)(399,0)
\path(50,54)(60,52)(70,50)(80,48)(90,46)(100,45)(110,46)(120,48)(130,50)(140,52)(150,54)
\path(70,50)(80,52)(90,54)(100,55)(110,54)(120,52)(130,50)
\path(300,54)(290,52)(280,50)(270,48)(260,46)(250,45)(240,46)(230,48)(220,50)(210,52)(200,54)
\path(280,50)(270,52)(260,54)(250,55)(240,54)(230,52)(220,50)
\path(350,54)(360,52)(370,50)(380,48)(390,46)(400,45)(410,46)(420,48)(430,50)(440,52)(450,54)
\path(370,50)(380,52)(390,54)(400,55)(410,54)(420,52)(430,50)
\put(250,22){\circle*{10}}
\put(50,115){$\psi^3-\lambda_1\psi^2+\lambda_2\psi-\lambda_3$}
\put(510,45){$=$}
\path(550,50)(551,57)(552,60)(555,66)(558,70)(564,76)(568,79)(571,81)(573,82)(576,84)(584,88)(587,89)(590,90)(602,94)(610,96)(615,97)(622,98)(630,99)(650,100)(651,100)
\path(550,50)(551,43)(552,40)(555,34)(558,30)(564,24)(568,21)(571,19)(573,18)(576,16)(584,12)(587,11)(590,10)(602,6)(610,4)(615,3)(622,2)(630,1)(650,0)(651,0)
\path(725,88)(710,90)(698,94)(690,96)(685,97)(678,98)(670,99)(650,100)(649,100)
\path(725,12)(710,10)(698,6)(690,4)(685,3)(678,2)(670,1)(650,0)(649,0)
\path(900,50)(899,57)(898,60)(895,66)(892,70)(886,76)(882,79)(879,81)(877,82)(874,84)(866,88)(863,89)(860,90)(848,94)(840,96)(835,97)(828,98)(820,99)(800,100)(799,100)
\path(900,50)(899,43)(898,40)(895,34)(892,30)(886,24)(882,21)(879,19)(877,18)(874,16)(866,12)(863,11)(860,10)(848,6)(840,4)(835,3)(828,2)(820,1)(800,0)(799,0)
\path(725,88)(740,90)(752,94)(760,96)(765,97)(772,98)(780,99)(800,100)(801,100)
\path(725,12)(740,10)(752,6)(760,4)(765,3)(772,2)(780,1)(800,0)(801,0)
\path(700,54)(690,52)(680,50)(670,48)(660,46)(650,45)(640,46)(630,48)(620,50)(610,52)(600,54)
\path(680,50)(670,52)(660,54)(650,55)(640,54)(630,52)(620,50)
\path(750,54)(760,52)(770,50)(780,48)(790,46)(800,45)(810,46)(820,48)(830,50)(840,52)(850,54)
\path(770,50)(780,52)(790,54)(800,55)(810,54)(820,52)(830,50)
\put(675,115){$\psi-\lambda_1$}
\put(950,50){\circle{100}}
\path(1000,50)(990,48)(980,46)(970,44)(960,42)(950,41)(940,42)(930,44)(920,46)(910,48)(900,50)
\path(990,52)(980,54)\path(970,56)(960,58)\path(950,59)(940,58)\path(930,56)(920,54)\path(910,52)(900,50)
\put(950,20){\circle*{10}}
\put(1100,50){\ellipse{200}{100}}
\path(1050,54)(1060,52)(1070,50)(1080,48)(1090,46)(1100,45)(1110,46)(1120,48)(1130,50)(1140,52)(1150,54)
\path(1070,50)(1080,52)(1090,54)(1100,55)(1110,54)(1120,52)(1130,50)
\put(1210,45){$+$}
\put(1350,0){\ellipse{200}{100}}
\path(1300,4)(1310,2)(1320,0)(1330,-2)(1340,-4)(1350,-5)(1360,-4)(1370,-2)(1380,0)(1390,2)(1400,4)
\path(1320,0)(1330,2)(1340,4)(1350,5)(1360,4)(1370,2)(1380,0)
\put(1500,0){\circle{100}}
\path(1550,0)(1540,-2)(1530,-4)(1520,-6)(1510,-8)(1500,-9)(1490,-8)(1480,-6)(1470,-4)(1460,-2)(1450,0)
\path(1540,2)(1530,4)\path(1520,6)(1510,8)\path(1500,9)(1490,8)\path(1480,6)(1470,4)\path(1460,2)(1450,0)
\put(1500,-30){\circle*{10}}
\put(1650,0){\ellipse{200}{100}}
\path(1600,4)(1610,2)(1620,0)(1630,-2)(1640,-4)(1650,-5)(1660,-4)(1670,-2)(1680,0)(1690,2)(1700,4)
\path(1620,0)(1630,2)(1640,4)(1650,5)(1660,4)(1670,2)(1680,0)
\put(1500,100){\ellipse{200}{100}}
\path(1450,104)(1460,102)(1470,100)(1480,98)(1490,96)(1500,95)(1510,96)(1520,98)(1530,100)(1540,102)(1550,104)
\path(1470,100)(1480,102)(1490,104)(1540,105)(1510,104)(1520,102)(1530,100)
\end{picture}
\end{eqnarray*}
}

Genus $4$: In $\ocM_{4,1}$,
{\tiny
\begin{eqnarray*}
\setlength{\unitlength}{0.0075cm}
\begin{picture}(1900,700)(-50,0)
\thicklines

\path(600,600)(601,607)(602,610)(605,616)(608,620)(614,626)(618,629)(621,631)(623,632)(626,634)(634,638)(637,639)(640,640)(652,644)(660,646)(665,647)(672,648)(680,649)(700,650)(701,650)
\path(600,600)(601,593)(602,590)(605,584)(608,580)(614,574)(618,571)(621,569)(623,568)(626,566)(634,562)(637,561)(640,560)(652,556)(660,554)(665,553)(672,552)(680,551)(700,550)(701,550)
\path(775,638)(760,640)(748,644)(740,646)(735,647)(728,648)(720,649)(700,650)(699,650)
\path(775,562)(760,560)(748,556)(740,554)(735,553)(728,552)(720,551)(700,550)(699,550)
\path(775,638)(790,640)(802,644)(810,646)(815,647)(822,648)(830,649)(850,650)(851,650)
\path(775,562)(790,560)(802,556)(810,554)(815,553)(822,552)(830,551)(850,550)(851,550)
\path(650,604)(660,602)(670,600)(680,598)(690,596)(700,595)(710,596)(720,598)(730,600)(740,602)(750,604)
\path(670,600)(680,602)(690,604)(700,605)(710,604)(720,602)(730,600)
\path(925,638)(910,640)(898,644)(890,646)(885,647)(878,648)(870,649)(850,650)(849,650)
\path(925,562)(910,560)(898,556)(890,554)(885,553)(878,552)(870,551)(850,550)(849,550)
\path(925,638)(940,640)(952,644)(960,646)(965,647)(972,648)(980,649)(1000,650)(1001,650)
\path(925,562)(940,560)(952,556)(960,554)(965,553)(972,552)(980,551)(1000,550)(1001,550)
\path(800,604)(810,602)(820,600)(830,598)(840,596)(850,595)(860,596)(870,598)(880,600)(890,602)(900,604)
\path(820,600)(830,602)(840,604)(850,605)(860,604)(870,602)(880,600)
\path(1075,638)(1060,640)(1048,644)(1040,646)(1035,647)(1028,648)(1020,649)(1000,650)(999,650)
\path(1075,562)(1060,560)(1048,556)(1040,554)(1035,553)(1028,552)(1020,551)(1000,550)(999,550)
\path(1075,638)(1090,640)(1102,644)(1110,646)(1115,647)(1122,648)(1130,649)(1150,650)(1151,650)
\path(1075,562)(1090,560)(1102,556)(1110,554)(1115,553)(1122,552)(1130,551)(1150,550)(1151,550)
\path(950,604)(960,602)(970,600)(980,598)(990,596)(1000,595)(1010,596)(1020,598)(1030,600)(1040,602)(1050,604)
\path(970,600)(980,602)(990,604)(1000,605)(1010,604)(1020,602)(1030,600)
\path(1250,600)(1249,607)(1248,610)(1245,616)(1242,620)(1236,626)(1232,629)(1229,631)(1227,632)(1224,634)(1216,638)(1213,639)(1210,640)(1198,644)(1190,646)(1185,647)(1178,648)(1170,649)(1150,650)(1149,650)
\path(1250,600)(1249,593)(1248,590)(1245,584)(1242,580)(1236,574)(1232,571)(1229,569)(1227,568)(1224,566)(1216,562)(1213,561)(1210,560)(1198,556)(1190,554)(1185,553)(1178,552)(1170,551)(1150,550)(1149,550)
\path(1100,604)(1110,602)(1120,600)(1130,598)(1140,596)(1150,595)(1160,596)(1170,598)(1180,600)(1190,602)(1200,604)
\path(1120,600)(1130,602)(1140,604)(1150,605)(1160,604)(1170,602)(1180,600)
\put(925,600){\circle*{10}}
\put(655,665){$\psi^4-\lambda_1\psi^3+\lambda_2\psi^2-\lambda_3\psi+\lambda_4$}

\put(1265,595){$=$}

\path(50,400)(51,407)(52,410)(55,416)(58,420)(64,426)(68,429)(71,431)(73,432)(76,434)(84,438)(87,439)(90,440)(102,444)(110,446)(115,447)(122,448)(130,449)(150,450)(151,450)
\path(50,400)(51,393)(52,390)(55,384)(58,380)(64,374)(68,371)(71,369)(73,368)(76,366)(84,362)(87,361)(90,360)(102,356)(110,354)(115,353)(122,352)(130,351)(150,350)(151,350)
\path(225,438)(210,440)(198,444)(190,446)(185,447)(178,448)(170,449)(150,450)(149,450)
\path(225,362)(210,360)(198,356)(190,354)(185,353)(178,352)(170,351)(150,350)(149,350)
\path(225,438)(240,440)(252,444)(260,446)(265,447)(272,448)(280,449)(300,450)(301,450)
\path(225,362)(240,360)(252,356)(260,354)(265,353)(272,352)(280,351)(300,350)(301,350)
\path(375,438)(360,440)(348,444)(340,446)(335,447)(328,448)(320,449)(300,450)(299,450)
\path(375,362)(360,360)(348,356)(340,354)(335,353)(328,352)(320,351)(300,350)(299,350)
\path(375,438)(390,440)(402,444)(410,446)(415,447)(422,448)(430,449)(450,450)(451,450)
\path(375,362)(390,360)(402,356)(410,354)(415,353)(422,352)(430,351)(450,350)(451,350)
\path(550,400)(549,407)(548,410)(545,416)(542,420)(536,426)(532,429)(529,431)(527,432)(524,434)(516,438)(513,439)(510,440)(498,444)(490,446)(485,447)(478,448)(470,449)(450,450)(449,450)
\path(550,400)(549,393)(548,390)(545,384)(542,380)(536,374)(532,371)(529,369)(527,368)(524,366)(516,362)(513,361)(510,360)(498,356)(490,354)(485,353)(478,352)(470,351)(450,350)(449,350)
\path(100,404)(110,402)(120,400)(130,398)(140,396)(150,395)(160,396)(170,398)(180,400)(190,402)(200,404)
\path(120,400)(130,402)(140,404)(150,405)(160,404)(170,402)(180,400)
\path(350,404)(340,402)(330,400)(320,398)(310,396)(300,395)(290,396)(280,398)(270,400)(260,402)(250,404)
\path(330,400)(320,402)(310,404)(300,405)(290,404)(280,402)(270,400)
\path(400,404)(410,402)(420,400)(430,398)(440,396)(450,395)(460,396)(470,398)(480,400)(490,402)(500,404)
\path(420,400)(430,402)(440,404)(450,405)(460,404)(470,402)(480,400)
\put(175,465){$\psi^2-\lambda_1\psi+\lambda_2$}
\put(600,400){\circle{100}}
\path(650,400)(640,398)(630,396)(620,394)(610,392)(600,391)(590,392)(580,394)(570,396)(560,398)(550,400)
\path(640,402)(630,404)
\path(620,406)(610,408)
\path(600,409)(590,408)
\path(580,406)(570,404)
\path(560,402)(550,400)
\put(600,370){\circle*{10}}
\put(750,400){\ellipse{200}{100}}
\path(700,404)(710,402)(720,400)(730,398)(740,396)(750,395)(760,396)(770,398)(780,400)(790,402)(800,404)
\path(720,400)(730,402)(740,404)(750,405)(760,404)(770,402)(780,400)

\put(860,395){$+$}

\path(900,400)(901,407)(902,410)(905,416)(908,420)(914,426)(918,429)(921,431)(923,432)(926,434)(934,438)(937,439)(940,440)(952,444)(960,446)(965,447)(972,448)(980,449)(1000,450)(1001,450)
\path(900,400)(901,393)(902,390)(905,384)(908,380)(914,374)(918,371)(921,369)(923,368)(926,366)(934,362)(937,361)(940,360)(952,356)(960,354)(965,353)(972,352)(980,351)(1000,350)(1001,350)
\path(1075,438)(1060,440)(1048,444)(1040,446)(1035,447)(1028,448)(1020,449)(1000,450)(999,450)
\path(1075,362)(1060,360)(1048,356)(1040,354)(1035,353)(1028,352)(1020,351)(1000,350)(999,350)
\path(1250,400)(1249,407)(1248,410)(1245,416)(1242,420)(1236,426)(1232,429)(1229,431)(1227,432)(1224,434)(1216,438)(1213,439)(1210,440)(1198,444)(1190,446)(1185,447)(1178,448)(1170,449)(1150,450)(1149,450)
\path(1250,400)(1249,393)(1248,390)(1245,384)(1242,380)(1236,374)(1232,371)(1229,369)(1227,368)(1224,366)(1216,362)(1213,361)(1210,360)(1198,356)(1190,354)(1185,353)(1178,352)(1170,351)(1150,350)(1149,350)
\path(1075,438)(1090,440)(1102,444)(1110,446)(1115,447)(1122,448)(1130,449)(1150,450)(1151,450)
\path(1075,362)(1090,360)(1102,356)(1110,354)(1115,353)(1122,352)(1130,351)(1150,350)(1151,350)
\path(950,404)(960,402)(970,400)(980,398)(990,396)(1000,395)(1010,396)(1020,398)(1030,400)(1040,402)(1050,404)
\path(970,400)(980,402)(990,404)(1000,405)(1010,404)(1020,402)(1030,400)
\path(1200,404)(1190,402)(1180,400)(1170,398)(1160,396)(1150,395)(1140,396)(1130,398)(1120,400)(1110,402)(1100,404)
\path(1180,400)(1170,402)(1160,404)(1150,405)(1140,404)(1130,402)(1120,400)
\put(1025,465){$\psi-\lambda_1$}
\put(1300,400){\circle{100}}
\path(1350,400)(1340,398)(1330,396)(1320,394)(1310,392)(1300,391)(1290,392)(1280,394)(1270,396)(1260,398)(1250,400)
\path(1340,402)(1330,404)
\path(1320,406)(1310,408)
\path(1300,409)(1290,408)
\path(1280,406)(1270,404)
\path(1260,402)(1250,400)
\put(1300,370){\circle*{10}}
\path(1350,400)(1351,407)(1352,410)(1355,416)(1358,420)(1364,426)(1368,429)(1371,431)(1373,432)(1376,434)(1384,438)(1387,439)(1390,440)(1402,444)(1410,446)(1415,447)(1422,448)(1430,449)(1450,450)(1451,450)
\path(1350,400)(1351,393)(1352,390)(1355,384)(1358,380)(1364,374)(1368,371)(1371,369)(1373,368)(1376,366)(1384,362)(1387,361)(1390,360)(1402,356)(1410,354)(1415,353)(1422,352)(1430,351)(1450,350)(1451,350)
\path(1525,438)(1510,440)(1498,444)(1490,446)(1485,447)(1478,448)(1470,449)(1450,450)(1449,450)
\path(1525,362)(1510,360)(1498,356)(1490,354)(1485,353)(1478,352)(1470,351)(1450,350)(1449,350)
\path(1700,400)(1699,407)(1698,410)(1695,416)(1692,420)(1686,426)(1682,429)(1679,431)(1677,432)(1674,434)(1666,438)(1663,439)(1660,440)(1648,444)(1640,446)(1635,447)(1628,448)(1620,449)(1600,450)(1599,450)
\path(1700,400)(1699,393)(1698,390)(1695,384)(1692,380)(1686,374)(1682,371)(1679,369)(1677,368)(1674,366)(1666,362)(1663,361)(1660,360)(1648,356)(1640,354)(1635,353)(1628,352)(1620,351)(1600,350)(1599,350)
\path(1525,438)(1540,440)(1552,444)(1560,446)(1565,447)(1572,448)(1580,449)(1600,450)(1601,450)
\path(1525,362)(1540,360)(1552,356)(1560,354)(1565,353)(1572,352)(1580,351)(1600,350)(1601,350)
\path(1400,404)(1410,402)(1420,400)(1430,398)(1440,396)(1450,395)(1460,396)(1470,398)(1480,400)(1490,402)(1500,404)
\path(1420,400)(1430,402)(1440,404)(1450,405)(1460,404)(1470,402)(1480,400)
\path(1650,404)(1640,402)(1630,400)(1620,398)(1610,396)(1600,395)(1590,396)(1580,398)(1570,400)(1560,402)(1550,404)
\path(1630,400)(1620,402)(1610,404)(1600,405)(1590,404)(1580,402)(1570,400)
\put(1475,465){$\psi-\lambda_1$}

\put(1710,395){$+$}

\put(-40,145){$+$}

\path(0,100)(1,107)(2,110)(5,116)(8,120)(14,126)(18,129)(21,131)(23,132)(26,134)(34,138)(37,139)(40,140)(52,144)(60,146)(65,147)(72,148)(80,149)(100,150)(101,150)
\path(0,100)(1,93)(2,90)(5,84)(8,80)(14,74)(18,71)(21,69)(23,68)(26,66)(34,62)(37,61)(40,60)(52,56)(60,54)(65,53)(72,52)(80,51)(100,50)(101,50)
\path(175,138)(160,140)(148,144)(140,146)(135,147)(128,148)(120,149)(100,150)(99,150)
\path(175,62)(160,60)(148,56)(140,54)(135,53)(128,52)(120,51)(100,50)(99,50)
\path(350,100)(349,107)(348,110)(345,116)(342,120)(336,126)(332,129)(329,131)(327,132)(324,134)(316,138)(313,139)(310,140)(298,144)(290,146)(285,147)(278,148)(270,149)(250,150)(249,150)
\path(350,100)(349,93)(348,90)(345,84)(342,80)(336,74)(332,71)(329,69)(327,68)(324,66)(316,62)(313,61)(310,60)(298,56)(290,54)(285,53)(278,52)(270,51)(250,50)(249,50)
\path(175,138)(190,140)(202,144)(210,146)(215,147)(222,148)(230,149)(250,150)(251,150)
\path(175,62)(190,60)(202,56)(210,54)(215,53)(222,52)(230,51)(250,50)(251,50)
\path(50,104)(60,102)(70,100)(80,98)(90,96)(100,95)(110,96)(120,98)(130,100)(140,102)(150,104)
\path(70,100)(80,102)(90,104)(100,105)(110,104)(120,102)(130,100)
\path(300,104)(290,102)(280,100)(270,98)(260,96)(250,95)(240,96)(230,98)(220,100)(210,102)(200,104)
\path(280,100)(270,102)(260,104)(250,105)(240,104)(230,102)(220,100)
\put(125,165){$\psi-\lambda_1$}
\put(400,100){\circle{100}}
\path(450,100)(440,98)(430,96)(420,94)(410,92)(400,91)(390,92)(380,94)(370,96)(360,98)(350,100)
\path(440,102)(430,104)
\path(420,106)(410,108)
\path(400,109)(390,108)
\path(380,106)(370,104)
\path(360,102)(350,100)
\put(400,70){\circle*{10}}
\put(550,100){\ellipse{200}{100}}
\path(500,104)(510,102)(520,100)(530,98)(540,96)(550,95)(560,96)(570,98)(580,100)(590,102)(600,104)
\path(520,100)(530,102)(540,104)(550,105)(560,104)(570,102)(580,100)
\put(400,200){\ellipse{200}{100}}
\path(350,204)(360,202)(370,200)(380,198)(390,196)(400,195)(410,196)(420,198)(430,200)(440,202)(450,204)
\path(370,200)(380,202)(390,204)(400,205)(410,204)(420,202)(430,200)

\put(660,145){$+$}

\put(800,150){\ellipse{200}{100}}
\path(750,154)(760,152)(770,150)(780,148)(790,146)(800,145)(810,146)(820,148)(830,150)(840,152)(850,154)
\path(770,150)(780,152)(790,154)(800,155)(810,154)(820,152)(830,150)
\put(950,150){\circle{100}}
\path(1000,150)(990,148)(980,146)(970,144)(960,142)(950,141)(940,142)(930,144)(920,146)(910,148)(900,150)
\path(990,152)(980,154)
\path(970,156)(960,158)
\path(950,159)(940,158)
\path(930,156)(920,154)
\path(910,152)(900,150)
\put(950,120){\circle*{10}}
\put(950,50){\ellipse{200}{100}}
\path(900,54)(910,52)(920,50)(930,48)(940,46)(950,45)(960,46)(970,48)(980,50)(990,52)(1000,54)
\path(920,50)(930,52)(940,54)(950,55)(960,54)(970,52)(980,50)
\put(950,250){\ellipse{200}{100}}
\path(900,254)(910,252)(920,250)(930,248)(940,246)(950,245)(960,246)(970,248)(980,250)(990,252)(1000,254)
\path(920,250)(930,252)(940,254)(950,255)(960,254)(970,252)(980,250)
\put(1100,150){\ellipse{200}{100}}
\path(1050,154)(1060,152)(1070,150)(1080,148)(1090,146)(1100,145)(1110,146)(1120,148)(1130,150)(1140,152)(1150,154)
\path(1070,150)(1080,152)(1090,154)(1100,155)(1110,154)(1120,152)(1130,150)

\put(1210,145){$-$}

\path(1250,150)(1251,157)(1252,160)(1255,166)(1258,170)(1264,176)(1268,179)(1271,181)(1273,182)(1276,184)(1284,188)(1287,189)(1290,190)(1302,194)(1310,196)(1315,197)(1322,198)(1330,199)(1350,200)(1351,200)
\path(1250,150)(1251,143)(1252,140)(1255,134)(1258,130)(1264,124)(1268,121)(1271,119)(1273,118)(1276,116)(1284,112)(1287,111)(1290,110)(1302,106)(1310,104)(1315,103)(1322,102)(1330,101)(1350,100)(1351,100)
\path(1425,188)(1410,190)(1398,194)(1390,196)(1385,197)(1378,198)(1370,199)(1350,200)(1349,200)
\path(1425,112)(1410,110)(1398,106)(1390,104)(1385,103)(1378,102)(1370,101)(1350,100)(1349,100)
\path(1600,150)(1599,157)(1598,160)(1595,166)(1592,170)(1586,176)(1582,179)(1579,181)(1577,182)(1574,184)(1566,188)(1563,189)(1560,190)(1548,194)(1540,196)(1535,197)(1528,198)(1520,199)(1500,200)(1499,200)
\path(1600,150)(1599,143)(1598,140)(1595,134)(1592,130)(1586,124)(1582,121)(1579,119)(1577,118)(1574,116)(1566,112)(1563,111)(1560,110)(1548,106)(1540,104)(1535,103)(1528,102)(1520,101)(1500,100)(1499,100)
\path(1425,188)(1440,190)(1452,194)(1460,196)(1465,197)(1472,198)(1480,199)(1500,200)(1501,200)
\path(1425,112)(1440,110)(1452,106)(1460,104)(1465,103)(1472,102)(1480,101)(1500,100)(1501,100)
\path(1300,154)(1310,152)(1320,150)(1330,148)(1340,146)(1350,145)(1360,146)(1370,148)(1380,150)(1390,152)(1400,154)
\path(1320,150)(1330,152)(1340,154)(1350,155)(1360,154)(1370,152)(1380,150)
\path(1550,154)(1540,152)(1530,150)(1520,148)(1510,146)(1500,145)(1490,146)(1480,148)(1470,150)(1460,152)(1450,154)
\path(1530,150)(1520,152)(1510,154)(1500,155)(1490,154)(1480,152)(1470,150)
\put(1650,150){\circle{100}}
\path(1700,150)(1690,148)(1680,146)(1670,144)(1660,142)(1650,141)(1640,142)(1630,144)(1620,146)(1610,148)(1600,150)
\path(1690,152)(1680,154)
\path(1670,156)(1660,158)
\path(1650,159)(1640,158)
\path(1630,156)(1620,154)
\path(1610,152)(1600,150)
\put(1650,120){\circle*{10}}
\put(1750,150){\circle{100}}
\path(1800,150)(1790,148)(1780,146)(1770,144)(1760,142)(1750,141)(1740,142)(1730,144)(1720,146)(1710,148)(1700,150)
\path(1790,152)(1780,154)
\path(1770,156)(1760,158)
\path(1750,159)(1740,158)
\path(1730,156)(1720,154)
\path(1710,152)(1700,150)
\put(1750,50){\ellipse{200}{100}}
\path(1700,54)(1710,52)(1720,50)(1730,48)(1740,46)(1750,45)(1760,46)(1770,48)(1780,50)(1790,52)(1800,54)
\path(1720,50)(1730,52)(1740,54)(1750,55)(1760,54)(1770,52)(1780,50)
\put(1750,250){\ellipse{200}{100}}
\path(1700,254)(1710,252)(1720,250)(1730,248)(1740,246)(1750,245)(1760,246)(1770,248)(1780,250)(1790,252)(1800,254)
\path(1720,250)(1730,252)(1740,254)(1750,255)(1760,254)(1770,252)(1780,250)

\end{picture}
\end{eqnarray*}
}

We also have calculated the formula for
$ \psi^5 - \lambda_1 \psi^4 + \lambda_2 \psi^3 - \lambda_3 \psi^2 +
\lambda_4 \psi - \lambda_5 $
in $\ocM_{5,1}$.
We do not write it here because it was calculated via a (possibly incorrect)
computer program and because it is rather long.
Note that non-integer coefficients do appear in genus $5$.

\end{document}